\documentclass[12pt,thmsa, reqno]{amsart}

\usepackage{mathrsfs}
\usepackage{amsfonts}
\usepackage{amssymb}
\usepackage{amsmath}

\usepackage[dvips]{graphics}
\input{amssym.def}
\input{amssym.tex}

\def\th{\theta}

\def\Cal{\mathcal}

\def\F{{\Cal F}}

\def\I{{\Cal I}}

\def\gnk{G_{n,k}}

\def\f0{f_0}
\def\Fc0{\varphi_0}
\def\rn{\bbr^n}

\def\I_k {I_{-}^{k/2}}
\def\I+k {I_{+}^{k/2}}

\def\bbr{{\Bbb R}}

\def\cos{{\hbox{\rm cos}}}

\def\gnk{G_{n,k}}

\def\part{\partial}
\def\intl{\int\limits}

\def\Gam{\Gamma}

\def\vp{\varphi}

\def\gam{\gamma}

\font\frak=eufm10

\def\fr#1{\hbox{\frak #1}}

\def\frM{\fr{M}}

\def\cos{{\hbox{\rm cos}}}

\def\part{\partial}
\def\intl{\int\limits}

\def\Gam{\Gamma}

\def\gnk{G_{n,k}}

\newtheorem{theorem}{Theorem}[section]

\newtheorem{corollary}[theorem]{Corollary}

\theoremstyle{remark}

\numberwithin{equation}{section}

\newcommand{\be}{\begin{equation}}
\newcommand{\ee}{\end{equation}}

\newcommand{\bea}{\begin{eqnarray}}
\newcommand{\eea}{\end{eqnarray}}
\newcommand{\Bea}{\begin{eqnarray*}}
\newcommand{\Eea}{\end{eqnarray*}}

\def\sideremark#1{\ifvmode\leavevmode\fi\vadjust{\vbox to0pt{\vss% the remark
 \hbox to 0pt{\hskip\hsize\hskip1em%                          will appear only
\vbox{\hsize2cm\tiny\raggedright\pretolerance10000%          on the side
 \noindent #1\hfill}\hss}\vbox to8pt{\vfil}\vss}}}%
\begin{document}

\title[On the Determination of Star Bodies]
{On the Determination of Star Bodies from Their Half-Sections}

\author{B. Rubin }
\address{Department of Mathematics, Louisiana State University, Baton Rouge,
Louisiana 70803, USA}
\email{borisr@math.lsu.edu}

%\thanks{ The work  was supported in part by
%the Edmund Landau Center for Research in Mathematical Analysis
%and Related Areas, sponsored by the Minerva Foundation (Germany).}


\date{April 16, 2016}

%\dedicatory{This paper is dedicated to our authors.}


\begin{abstract}

We obtain explicit inversion formulas for the Radon-like transform that assigns to a function on the unit sphere   the integrals of that  function  over    hemispheres lying in lower dimensional central cross-sections. The results are applied to determination of  star bodies from the volumes of their central half-sections.
 \end{abstract}

\maketitle

\section{Introduction}
\setcounter{equation}{0}

Let $K$ be a   compact subset of $\rn$ which contains the origin  $o$ as an  interior point and is star-shaped with respect to  $o$. Such a set $K$ is called a
 {\it star body} if the  radial function
\be\label {rad}\rho_K(\th)=\max \{r\ge 0 : r\th \in K\}, \qquad \th\in S^{n-1},\ee
that determines the shape of  $K$ is continuous; see Gardner \cite{Ga06}.
 We denote by $\gnk$,  $2\le k\le n-1$,  the Grassmann manifold of all $k$-dimensional  linear subspaces $\xi$ of $\rn$. Passing to spherical coordinates, one can evaluate the $k$-dimensional volume of  $K\cap \xi$ as
\be\label {vol} vol_k (K\cap \xi)=\frac{1}{k} \intl_{S^{n-1} \cap \xi} \rho_K^k(\th)\, d_\xi \th,\ee
 where $d_\xi \th$ stands for the corresponding $(k-1)$-dimensional surface element. The right-hand side of (\ref{vol}) is a constant multiple of the  Funk-type transform
 \be\label {kFunk}
(F_k f)(\xi)=\intl_{S^{n-1}\cap \xi} f(\th)\, d_\xi \th, \qquad \xi \in \gnk,
\ee
that can be explicitly inverted by different ways provided that $f$ is an even integrable function on $S^{n-1}$; see, e.g., \cite{He11, Ru02, Ru15b}. This fact makes it possible to determine the shape of $K$ from the knowledge of the volumes of $K\cap \xi$ for all $\xi\in \gnk$ provided that $K$ is origin-symmetric.
The star bodies which are not origin-symmetric cannot be uniquely reconstructed from the integrals (\ref{vol}); see, e.g., Gardner \cite[Section 7.2]{Ga06}. This statement agrees with a known fact that the  kernel of  $F_k$  on $L^1(S^{n-1})$ consists of odd functions.

In the case $k=n-1$, Groemer \cite{Gr98} who followed some ideas from  Backus \cite {B},
suggested to replace the hyperplane central sections  by the ``half-sections" $K\cap H(u,v)$, where
\be\label {vol1} H(u,v)=\{x \in \rn : \; x\cdot u=0,\; x\cdot v \ge 0\}, \qquad u\perp v,\ee
 is a half-plane determined by  mutually orthogonal unit vectors $u$ and $v$, so that the origin $o$ lies on the boundary of  $H(u,v)$.
 Thus (\ref{vol}) is substituted by  the lower dimensional hemispherical integral
  \be\label {vol5} vol_{n-1} (K\cap H(u,v))=\frac{1}{n-1} \intl_{S^{n-1} \cap H(u,v)} \rho_K^{n-1}(\th)\, d_\xi \th.\ee
One of the main results of \cite{Gr98} states that if  the star bodies $K$ and $L$  satisfy
\[ vol_{n-1} (K\cap H(u,v))=vol_{n-1} (L\cap H(u,v)) \]
 for all $u,v\in  S^{n-1}$, $u\perp v$,
then $K=L$.  This uniqueness result was extended by Goodey and Weil \cite{GW06} to half-sections of arbitrary dimension $2\le k\le n-1$.

To the best of our knowledge, the following important questions remained open.

{\bf Question 1.} {\it How can we explicitly reconstruct the shape of  $K$ (or the radial function $\rho_K$) from volumes (\ref{vol5}) or, more generally, from the corresponding  $k$-dimensional volumes?}

{\bf Question 2.}  {\it  How can we  eliminate overdeterminedness of the inversion problem in Question 1?}

Regarding Question 2, we observe that  (\ref{vol1}) parameterizes the  corresponding half-sections by the elements of the Stiefel manifold $V_{n,2}=\{(u,v): u,v\in  S^{n-1}, u\perp v\}$, so that $\dim V_{n,2}=2n-3>n-1=\dim S^{n-1}$ if $n>2$. In other words, the dimension of the target space is greater than the dimension of the source space. The latter means that the  inversion problem in Question 1 can  be overdetermined. In the case $k<n-1$, the difference between the dimensions of the target space and  source space is  bigger because the corresponding lower dimensional Funk transform (\ref{kFunk}) is overdetermined itself. This situation is pursuant to
  Gel'fand's celebrated question \cite{Gelf60} on  how to reduce  overdeterminedness  of transformations in integral geometry. In our case it means the following

\noindent  {\bf Problem.} {\it Find an $(n-1)$-dimensional submanifold $\frM$ of the manifold of all $k$-dimensional  central half-sections of $K$ so that  $K$
  could be recovered from the volumes of the half-sections belonging to $\frM$ only.}

In the present article we solve this problem give and an answer to Questions 1 and 2. The basic idea is to  consider half-sections lying in the open  half-spaces
\be\label{hpm} H_{\pm}=\{x =(x_1, \ldots, x_n)\in\rn : \; \pm \,x_n >0\}\ee
 separately. The case of all  $2\le k\le n-1$ is considered in Section 2. Here our inversion formulas remain overdetermined if $k<n-1$. This overdeterminedness is eliminated in Section 3.

{\bf Acknowledgements.} The author is grateful to professors Richard Gardner,  Paul Goodey, and Wolfgang Weil for useful discussion.

\section{Inversion formulas}

We consider the  following hemispherical modifications   of the Funk-type transform (\ref {kFunk}) when a function $f\in L^1(S^{n-1})$ is integrated over $(k-1)$-dimensional hemispheres $S^{n-1}_{\pm}\cap \xi$, $\xi \in \gnk$, lying in the $(n-1)$-dimensional hemispheres
\[
S^{n-1}_{\pm}=\{\th =(\th_1, \ldots, \th_n)\in S^{n-1} : \; \pm\th_n >0\},\]
respectively. Specifically, for any subspace $\xi \in \gnk$, not orthogonal to $e_n=(0, \ldots, 0,1)$, we set
\be\label {kFunk+}  (F_k^+ f)(\xi)=\intl_{S^{n-1}_+\cap \xi} f(\th)\, d_\xi \th \qquad \left (=\intl_{S^{n-1}_-\cap \xi} f(-\th)\, d_\xi \th\right ),\ee
\be\label {kFunk-}  (F_k^- f)(\xi)=\intl_{S^{n-1}_-\cap \xi} f(\th)\, d_\xi \th \qquad \left (=\intl_{S^{n-1}_+\cap \xi} f(-\th)\, d_\xi \th\right ).\ee
Clearly,
\be\label {kFunk+1}
(F_k^+ f)(\xi)=\frac{1}{2}(F_k f_1)(\xi), \qquad f_1 (\th)=\left \{\begin{array} {ll} f(\th) &\mbox{if $\th\in S^{n-1}_+$,}\\
f(-\th) &\mbox{if $\th\in S^{n-1}_-$,}\end{array}\right.\ee
\be\label {kFunk-1}
(F_k^- f)(\xi)=\frac{1}{2}(F_k f_2)(\xi), \qquad f_2 (\th)=\left \{\begin{array} {ll} f(-\th) &\mbox{if $\th\in S^{n-1}_+$,}\\
f(\th) &\mbox{if $\th\in S^{n-1}_-$,}\end{array}\right.\ee
where  $f_i$ $(i=1,2)$ are integrable even functions on  $S^{n-1}$ that can be  reconstructed from $\vp^{\pm} =F_k^{\pm} f$ by the formulas $f_1=2F_k^{-1} \vp^+$ and $f_2=2F_k^{-1} \vp^-$. Combining these formulas, we obtain
\be\label {kFunk1}
f(\th)=\left \{\begin{array} {ll} 2(F_k^{-1} \vp^+)(\th) &\mbox{if $\th\in S^{n-1}_+$,}\\
2(F_k^{-1} \vp^-)(\th) &\mbox{if $\th\in S^{n-1}_-$.}\end{array}\right.\ee
Thus we have proved the following
\begin{theorem} \label{theorem1} For $2\le k\le n-1$, a function $f\in L^1(S^{n-1})$ can be recovered from the  integrals
$\vp^{\pm} =F_k^{\pm} f$ by the formula (\ref{kFunk1}).
\end{theorem}

Some comments are in order.

{\bf 1.} Let $\xi \in \gnk $ be a subspace which is not orthogonal to $e_n=(0, \ldots, 0,1)$. We denote $\xi_{\pm}=\xi \cap \{x\in \rn : \pm \,x_n >0\}$. Then Theorem \ref{theorem1} implies the following

\begin{corollary} \label{corollary1}  The radial function $\rho_K$ of the star body $K$ can be recovered from the volumes $v^{\pm} (\xi)=vol_k (K \cap \xi_{\pm})$ by the formula
 \be\label {kFunk1b}
\rho_K^k(\th)=\left \{\begin{array} {ll} 2k(F_k^{-1} v^+)(\th) &\mbox{if $\th\in S^{n-1}_+$,}\\
2k(F_k^{-1} v^-)(\th) &\mbox{if $\th\in S^{n-1}_-$.}\end{array}\right.\ee
For $\th_n=0$ it can be determined from (\ref{kFunk1b}) by continuity.
\end{corollary}

{\bf 2.} In the case $k=n-1$, we can set $\xi=u^\perp$, $u \in S^{n-1}$, and write   (\ref{kFunk}) as the usual Funk transform
\be\label {Funk}
(Ff)(u)=\intl_{\{\th \in S^{n-1}:\; u\cdot \th=0\}} f(\th)\, d_u \th.
\ee
Similar notations $(F^{\pm} f)(u)$ can be used for the hemispherical transforms (\ref{kFunk+}) and (\ref{kFunk-}).

\vskip 0.2 truecm

{\bf 3.}  The above reasoning shows that to reconstruct $f$ on $S^{n-1}_+$ (or $S^{n-1}_-$) the knowledge of $\vp^{+} =F_k^{+} f$ (or $\vp^{-} =F_k^{-} f$, resp.) is sufficient.

\vskip 0.2 truecm

{\bf 4.}  A theory of the Funk transform  provides a variety of inversion formulas for $F_k$; see, e.g., \cite {GGG, He11, Ru02, Ru13b, Ru15b}. The functions $f_1$ and $f_2$ in (\ref{kFunk+1}) and (\ref{kFunk-1}) can have  discontinuity on the equator $\th_n=0$. It means that we cannot apply inversion formulas for $F_k$ (at least, straightforward) in which the smoothness is crucial. However, if $f\in L^p(S^{n-1})$, $1\le p< \infty$, then  $f_1$ and $f_2$ also belong to $L^p(S^{n-1})$ and can be reconstructed, e.g., by   the method of mean value operators  as follows.

 For $r=\cos\, \psi$, $\psi \in [0, \pi/2]$, and $\th \in S^{n-1}$, consider the shifted dual
 transform
\be\label {ufus}
(F^*_{k,\th} \vp )(r)= \intl_{d(\th, \xi) = \psi} \varphi(\xi) \,
d\mu(\xi),\ee
where $\vp$ is a function on $\gnk$, $d(\cdot, \cdot)$ denotes the geodesic distance, and $d\mu(\xi)$ stands for the corresponding canonical measure; see  \cite [Section 5] {Ru13b} for details.

\begin{theorem}\label{invrhys} \cite [Theorem 5.3] {Ru13b}
An even  function $f \in L^p (S^{n-1})$, $1\le p<\infty$,  can be recovered from $\vp=F_k f$ by the formula
\be\label{90ashel}
f(\th) \!=  \! \lim\limits_{s\to 1}  \left (\frac {1}{2s}\,\frac {\partial}{\partial s}\right )^k \left [\frac{\pi^{-k/2}}{\Gam (k/2)}\intl_0^s
(s^2 \!- \!r^2)^{k/2-1} \,(F^*_{k,\th} \vp)(r) \,r^k\,dr\right ].\ee
In particular, for $k$ even,
\be\label{90ashele}
f(\th) \!=  \! \lim\limits_{s\to 1} \frac{1}{2\pi^{k/2}} \left (\frac {1}{2s}\,\frac {\partial}{\partial s}\right )^{k/2}[s^{k-1}(F^*_{k,\th} \vp) (s)].\ee
Altenatively,
\be\label{90ashys}
f(\th) \!=  \! \lim\limits_{s\to 1} \, \left (\frac {\partial}{\partial s}\right )^k \left [\frac{2^{-k}\, \pi^{-k/2}}{\Gam (k/2)}\intl_0^s
(s^2 \!- \!r^2)^{k/2-1} (F^*_{k,\th} \vp)(r) \,dr\right ].\ee
The  limit in these formulas is understood in the $L^p$-norm.
\end {theorem}

\section{The Case $k<n-1$}

If $k<n-1$, the inversion problem for  $F_k^{\pm}$ and $F_k$ is overdetermined because the dimension of the target space is greater than  the dimension of the source space:
 \[\dim \gnk =k(n-k)> n-1=\dim S^{n-1}.\]  Below we eliminate this overdeterminedness by choosing an $(n-1)$-dimensional submanifold $\frM$ of $\gnk$ which is sufficient to reconstruct $f$ from  $(F_k^{\pm} f)(\xi)$ with  $\xi \in \frM$.  We proceed as in \cite{Ru15a}.

 Suppose that $\{e_1, \ldots, e_n\}$ is a standard orthonormal basis in $\rn$ and denote
\be\label {ttt609bu0cc} \bbr^{n-k}=\bbr e_1 \oplus \cdots \oplus  \bbr e_{n-k}, \qquad \bbr^{k}=\bbr e_{n-k+1}\oplus \cdots \oplus  \bbr e_{n},\ee
\be\label {hoo34d4pp}
\bbr^{k+1}=\bbr e_{n-k} \oplus \bbr^{k},  \quad  S^{k}=S^{n-1} \cap \bbr^{k+1}, \quad  S_{\pm}^{k}= S_{\pm}^{n-1} \cap \bbr^{k+1}.   \ee
Given a point $ v\in S^{n-k-1}=S^{n-1} \cap \bbr^{n-k}$, we fix
 an  orthogonal transformation $\gam_v$ in   $\bbr^{n-k}$,  so that $\gam_v  e_{n-k}=v$. Let
\be\label {hfos4609bu5a}
\tilde \gam_v =\left[\begin{array}{ll}  \gam_v &0\\
0& I_{k}
\end{array}\right],\ee
where $I_k$ is the identity $k\times k$ matrix.
We denote by $G_k (\tilde \gam_v \bbr^{k+1})$ the Grassmannian of all $k$-dimensional linear subspaces of  $\tilde \gam_v \bbr^{k+1}$ and set
\be\label {ttOO} \frM=\bigcup\limits_{v\in S^{n-k-1}} G_k (\tilde \gam_v \bbr^{k+1}).\ee
The restrictions $\tilde F_k^{\pm} f$ of $F_k^{\pm} f$ onto $\frM$ can be identified with  functions on the set
$\tilde S_{n,k}=\{(v,w): \, v\in S^{n-k-1}, \; w\in  S^{k}\}$. These functions are represented by the hemispherical integrals
\be\label {h9844pp}(\tilde F_k^{\pm} f)(v,w)=\intl_{ \{\eta\in S_{\pm}^{k}:\, \eta \cdot w=0\}}   f_v(\eta)\, d_w \eta, \qquad   f_v(\eta)=f (\tilde \gam_v \eta),\ee
 over $(k-1)$-dimensional hemispheres on  $S^{k}$. Thus  $\tilde F_k^{\pm}$ can be inverted as in the previous section.

Let us explain the details. We equip $\tilde S_{n,k}$ with the product measure $dv dw$, where
 $dv$ and $dw$ stand for the corresponding  surface elements on $S^{n-k-1}$ and $S^{k}$.  Using \cite[Theorem 3.2]{Ru15a} and taking into account that the restrictions of $f_v$ onto $S^k_{\pm}$ and their even extensions (cf. (\ref{kFunk+1}), (\ref{kFunk-1})) belong to the same $L^p$ spaces, we obtain the following existence result.

  \begin{theorem}  Let   $1\le k< n-1$, $f\in L^p(S^n)$, $n-k<p\le \infty$. Then $(\tilde F_k^{\pm} f)(v,w)$ are finite for almost all $(v,w)\in \tilde S_{n,k}$. If $p\le n-k$, then there are functions  $\tilde f_\pm\in L^p(S^{n-1})$ for which $(\tilde F_k^{\pm} f_\pm)(v,w)=\infty$.
\end{theorem}

Owing to  (\ref{kFunk1}), to reconstruct $f$  from
$\vp_v^\pm=\tilde \F_k^{\pm}f$,  it suffices to invert the usual co-dimension one Funk transform $\tilde F$ on $S^k$. This gives
\be\label {kFunk1v}
f_v(\eta)=\left \{\begin{array} {ll} 2(\tilde F^{-1} \vp_v^+)(\eta) &\mbox{if $\eta\in S^{k}_+$,}\\
2(\tilde F^{-1} \vp_v^-)(\eta) &\mbox{if $\eta\in S^{k}_-$.}\end{array}\right.\ee

If $f$ is a continuous function, its value
 at a point $\th \in S^{n-1}$ can be found as follows. We interpret  $\th$ as a column vector
$\th=(\th_1, \ldots, \th_{n+1})^T$ and set
\[\th'=(\th_1, \ldots, \th_{n-k})^T\in \bbr^{n-k},\qquad \th''=(\th_{n-k+1}, \ldots, \th_{n})^T\in \bbr^{k},\]
Suppose that $\th'\neq 0$ and set
\be\label {YYu6} v=\th'/|\th'|\in S^{n-k-1}, \quad \eta=(0, \ldots, 0,|\th'|, \th'')^T\in S^{k}.\ee
Then $\tilde \gam_v \eta=\th$ and we get
$f(\theta)= (\F^{-1} \vp_v)(\eta)$ for $v$ and $\eta$ as in (\ref{YYu6}).
If $\th'=0$, then $f(\theta)$ can be reconstructed by continuity from its values at the neighboring points.

 If $f$ is an arbitrary  function in $ L^p(S^n)$,  $n-k<p\le \infty$, then  $f$ can be explicitly reconstructed at almost all points on almost all spheres $S^{k}_v=\tilde \gam_v S^{k}$ by making use of known inversion formulas for the Funk transform on this class of functions; see, e.g., \cite[p. 297]{Ru15a}.

 The above reasoning is obviously applicable to reconstruction of the radial function $\rho_K$ of the star body $K$ from the volumes $v^{\pm} (\xi)=vol_k (K \cap \xi_{\pm})$ with $\xi \in \frM$, as in Corollary \ref{corollary1}.


\begin{thebibliography}{[ASMR]}

\bibitem {B} G. Backus,  Geographical  interpretation  of  measurements  of  average  phase velocities of surface waves over great circular and semi circular paths. Bull. Seism. Soc. Amer., 1964,  \textbf{54}, 571-610.


\bibitem {Ga06} R. J. Gardner,  \textit{Geometric Tomography (second edition)}. Cambridge University Press, New York, 2006.


\bibitem  {Gelf60}   I. M. Gel'fand,   \textit{Integral geometry and its relation to the theory of representations}.  Russian Math. Surveys \textbf{15} (1960), no. 2, 143--151.

\bibitem  {GGG} I. M. Gelfand, S. G. Gindikin, and  M. I. Graev, \textit{Selected Topics in Integral Geometry}, Translations of Mathematical Monographs, AMS,
Providence, Rhode Island, 2003.

\bibitem {GW06} P. Goodey and W. Weil, Average section functions for star-shaped sets. Adv. in Appl. Math. \textbf{36} (2006), no. 1, 70--84.

\bibitem {Gr98} H. Groemer, On a spherical integral transformation and sections of star bodies. Monatsh. Math. 126 (1998), no. 2, 117--124.

\bibitem  {He11}  S. Helgason,     Integral geometry and Radon transform,
Springer, New York-Dordrecht-Heidelberg-London, 2011.

\bibitem  {Ru02} B. Rubin, Inversion formulas for the spherical  Radon transform and the generalized cosine transform, Advances in Appl. Math.  \textbf{29} (2002), 471--497.

\bibitem  {Ru13b} \bysame,  On the Funk-Radon-Helgason inversion method in integral geometry, Contemp. Math.  \textbf{599} (2013), 175--198.

\bibitem  {Ru15a} \bysame,  Overdetermined transforms in integral geometry. Complex analysis and dynamical systems VI. Part 1, 291--313, Contemp. Math., \textbf{653}, Amer. Math. Soc., Providence, RI, 2015.

\bibitem  {Ru15b}  \bysame,  Introduction to Radon transforms (with elements of fractional calculus and harmonic analysis),  Cambridge University Press,  New York, 2015.



\end{thebibliography}
\end{document}